\documentclass[a4paper,11pt]{article}
\usepackage[frenchb,english]{babel}
\usepackage[latin1]{inputenc}
\usepackage{ifthen}
\usepackage{ifpdf}
\usepackage{authblk}
\usepackage[latin1]{inputenc}
\usepackage{url}

\ifpdf
\usepackage{aeguill}
\usepackage[pdftex]{graphicx,color}
\usepackage[pdfauthor={Dai Pra, Louis, Minelli},
            pdftitle={Synchronization via interacting reinforcement}]{hyperref}
\else
\usepackage[T1]{fontenc}
\usepackage[dvips]{graphicx}
\fi

\usepackage{dsfont}
\usepackage{eurosym}
\usepackage{latexsym}
\usepackage{amsmath,amssymb,amscd,amsfonts}
\usepackage{mathrsfs}

\usepackage{amsthm}
\usepackage{amstext}
\usepackage{xspace}
\usepackage{url}
\usepackage{enumerate}
\usepackage[all,2cell,ps]{xy}
\usepackage[outerbars]{changebar}
\usepackage{framed}

\theoremstyle{definition}
\newtheorem{Def}{Definition}[section]

\theoremstyle{plain}
\newtheorem{theorem}{Theorem}

\theoremstyle{plain}
\newtheorem{prop}{Proposition}[section]

\theoremstyle{plain}
\newtheorem{lemma}{Lemma}

\theoremstyle{plain}

\theoremstyle{remark}
\newtheorem{remark}{Remark}

\theoremstyle{remark}

\newcommand{\E}{\mathds{E}}
\newcommand{\Prob}{\mathds{P}}

\newcommand{\N}{\mathbb{N}}

\newcommand{\Var}{\mathrm{Var}}

\newcommand{\bd}[1]{\begin{Def}\label{#1}}
\newcommand{\ed}{\end{Def}}
\newcommand{\bp}[1]{\begin{prop}\label{#1}}
\newcommand{\ep}{\end{prop}}
\newcommand{\bt}[1]{\begin{theorem}\label{#1}}
\newcommand{\et}{\end{theorem}}
\newcommand{\be}[1]{\begin{equation}\label{#1}}
\newcommand{\ee}{\end{equation}}
\newcommand{\epr}{\end{proof}}
\newcommand{\bpr}{\begin{proof}}
\newcommand{\bl}[1]{\begin{lemma}\label{#1}}
\newcommand{\el}{\end{lemma}}
\newcommand{\eex}{\end{example}}
\newcommand{\br}[1]{\begin{remark}\label{#1}}
\newcommand{\er}{\end{remark}}

\def\ra{\rightarrow}

\def\a{\alpha}
\def\ra{\rightarrow}

\author[1]{P. Dai Pra\thanks{daipra@math.unipd.it}}
\author[2]{P.-Y. Louis\thanks{pierre-yves.louis@math.cnrs.fr}}
\author[3]{I.~G. Minelli\thanks{ida.minelli@dm.univaq.it}}
\affil[1]{Dipartimento di Matematica Pura ed Applicata,
	Universit\`a degli Studi di Padova,
	via Trieste 63,
	I-35121 Padova,
	Italy}
\affil[2]{Laboratoire de Math\'ematiques et Applications, UMR 7348,  Universit\'e de Poitiers|CNRS, 11 Boulevard Marie et Pierre Curie\\
86962 Technopole du Futuroscope de Poitiers -- Chasseneuil Cedex
France}
\affil[3]{Dipartimento di Ingegneria e Scienze dell'Informazione e Matematica, Universit\`a degli Studi dell'Aquila,
Via Vetoio (Coppito 1), 67100 Coppito (AQ), Italy}

\title{Synchronization via interacting reinforcement}
\date{Post-print of J. Appl. Prob. 51, 556--568 (2014) \\
doi:10.1239/jap/1402578643 \\
Copyright Applied Probability Trust 2014 \\
Creative Commons Attribution Non-Commercial No Derivatives License \\
received 10.01.2013, revision received 8.7.2013}

\begin{document}
\maketitle

\url{https://projecteuclid.org/euclid.jap/1402578643#info+}

\begin{abstract}
We consider a system of urns of Polya-type, with balls of two colors; the reinforcement of each urn depends both on the content of the same urn and on the average content of all urns. We show that the urns synchronize almost surely, in the sense that the fraction of balls of a given color converges almost surely, as the time goes to infinity, to the same limit for all urns. A normal approximation for a large number of urns is also obtained.

\end{abstract}

\section{Introduction}

Synchronization phenomena for stochastic systems have raised a considerable interest in the last decade, both for their impact in applications and their intrinsic mathematical beauty. We make no attempt of giving an overview of the many lines of research related to synchronization; we rather mention only two of them, which are related but developed mostly independently.

The word {\em synchronization} is often referred to {\em phase synchronization} of pulsating systems. Interesting examples are provided by {\em cellular dynamics}. Concentration of certain molecules within cells vary over time due to chemical reaction, and may lead to a periodic behavior. Chemical reactions may also involve molecules in the extracellular space, producing interaction between different cells: phase synchronization is actually observed when the interaction is sufficiently strong (see \cite{prindle2011sensing} for an account and a remarkable application). For this and other applications, a stochastic model that reasonably approximates the dynamics of phases is a noisy version of the {\em Kuramoto} model (\cite{kuramoto2003chemical, acebron2005kuramoto}). It is a system of $N$ coupled random rotators that, for a sufficiently small coupling parameter show decoherence of phases: if $\theta_j(t)$ denotes the phase of the $j$-th oscillator at time $t$, then the distribution of
\[
\rho := \left| \frac1N \sum_{j=1}^N e^{i \theta_j(t)} \right|
\]
is, for large $N$ and large $t$, concentrated near zero. Such decoherence is lost as the coupling parameter is increased beyond a critical value.

In other contexts, the notion of synchronization is used even in absence of an underlying periodic motion. In some models for metastable behavior (see e.g. \cite{Gentz, Bovier}), one considers a system of many particles subject to a symmetric, double-well potential, and interacting through an attractive potential that, in the models that are more mathematically treatable, is assumed of {\em mean field} type. At sufficiently low temperature, the majority of particles are found in the same well, thus breaking the symmetry of the model. Synchronization, in this context, is this coherent behavior of the majority.

In this paper we also consider a system whose components are subject to an attractive, mean field interaction. We however add an ingredient to the model: {\em reinforcement}. Generally speaking, by reinforcement in a stochastic dynamic we mean any mechanism for which the probability that a given transition occurs has an increasing dependence on the number of times that ``similar'' transitions has taken place in the past. The most elementary model of reinforcement is the {\em Polya urn}.

The simplest Polya urn model consists of an urn which contains balls of two different colors (for example, {at time $t=0$}, $b$ white and $a$ red balls). At each discrete time a ball is drawn out and it is replaced in the urn together with another ball of the same color.\\
Let $Z_t$ be the fraction of red balls at time $t$, namely, the conditional probability of drawing a red ball at time $t+1$, given the fraction of red balls at time $t$.\\
A well known result (see for instance~\cite{klenke2007probability}) says that $Z_t$ is a bounded martingale and in particular $\lim_{t\rightarrow \infty}Z_t=Z_\infty$ a.s., where $Z_\infty$ has Beta distribution with parameters $a$ and $b$.\\

Here we present a modified model, in which we consider a set of $N$ Polya's urns and we introduce a ``mean field interaction'' among them:
\begin{itemize}
\item
at time $0$, each urn contains $a$ red and $b$ white balls, $a \geq 1$, $b \geq 1$;
\item
at each time $t+1$, given the fraction $Z_t(j)$, for $j=1,\ldots,N$, of red balls in each urn $j$ at time $t$, independently of what happens in all the other urns, a new red ball is replaced in urn $i$ with conditional probability
$\alpha Z_t+(1-\alpha)Z_t(i)$, where
$Z_t$ is the total fraction of red balls at time $t$, i.e. $Z_t=\frac{1}{N}\sum_{i=1}^N Z_t(i)$ and $\alpha\in [0,1]$.
\end{itemize}

The case $\alpha=0$ corresponds to $N$ independent copies of the classical Polya's urn described above. Thus, for $i=1,2,\ldots,N$, the $Z_t(i)$'s converge, as $t \rightarrow +\infty$, to i.i.d. random variables, whose distribution is Beta with parameters $a,\ b$. We show that, for $\alpha >0$, $Z_i(t)$ is no longer a martingale; however its almost sure limit as $t \ra +\infty$ exists, and {\em it is the same for all} $i=1,2,\ldots,N$. We refer to this phenomenon as {\em almost sure synchronization}. We also obtain bounds on the $L^2$ distance between $Z_i(t)$ and $Z_j(t)$, $i \neq j$, which goes to zero {\em uniformly} in $N$ and as a power law in $t$.

We remark that more general and complex models of interacting Polya urns have been considered recently by various authors (see e.g. \cite{pemantle2007survey, cirillo2012polya, paganoni2004interacting, marsili1998self}). In particular, urns with mean-field interaction have been considered in \cite{launay2011interacting,launay2012generalized}, but with a different  reinforcement scheme: the probability of drawing a red (or black) ball is proportional to the exponential of the number of red balls, rather than to the number of red balls, leading to a quite different synchronization picture.

The paper is organized as follows. In Section \ref{sec:model} we give some basic properties of the model, while Section \ref{sec:sync} is devoted to synchronization. Finally, in Section \ref{sec:CLT} we prove a Central Limit Theorem for $Z(t)$, in the limit as $N \ra +\infty$.

\section{The model} \label{sec:model}
\subsection{Definition} \label{subsec:def}
Consider a probability space in which a family $\{ U(t,i);\  t,i\in \mathbb N \}$
of i.i.d. random variables with uniform distribution on $[0,1]$ is defined.\\
For example, we take
$\Omega=[0,1]^{\N^2},\ \ \mathcal{F}=\mathcal{B}([0,1])^{\otimes \N^2},\ \ \mathbb{P}=\lambda^{\otimes \N^2}$ (where $\mathcal{B}([0,1])$ and $\lambda$ denote respectively the Borel sigma-algebra and the Lebesgue measure on $[0,1]$), and the $U(t,i)$'s are the canonical projections.
We pose $\mathcal{F}_t=\sigma\big( U(s,i);\  0\leq s \leq t, i \in \N \big)$.\\

Let $m=a+b$ be the total number of balls in a single urn at time $0$.\\
In what follows we shall denote respectively by $X^N_t(i)$ and $Z^N_t(i)=\frac{X^N_t(i)}{t+m}$ the number and the fraction of red balls in urn $i$ at time $t$ for the model with $N$ urns.\\
Then, $X^N_t=\sum_{i=1}^N X^N_t(i)$ will denote the total number of red balls at time $t$ and
$Z^N_t=\frac{1}{N}\sum_{i=1}^N Z^N_t(i)=\frac{X_t}{N(t+m)}$ will be the total fraction of red balls at time $t$ in the model with $N$ urns.\\

Now, for a fixed $N\in \N$ we define for each $i=1,\ldots,N$
\begin{eqnarray*}
&&X_0^N(i)=1,\ \ \ Z_0^N(i)=\frac{a}{m}\\
&&X_{t+1}^N(i)=X^N_t(i)+ Y_t(i),\ \ \ Z^N_{t+1}(i)=\frac{t+m}{t+m+1}Z^N_{t}(i)+\frac{1}{t+m+1}Y_t(i)\\
\end{eqnarray*}

where $Y_t(i)=I_{\{U(t+1,i)\leq \alpha Z^N_t+(1-\alpha)Z^N_t(i)\}}$.\\

A simple calculation shows that the random variables $Y_t(i)$ for $i=1,\ldots,N$ are conditionally independent given $\mathcal{F}_t$ and have conditional distribution which is Bernoulli with parameter $\alpha Z_t+(1-\alpha)Z_t(i)$.\\

\subsection{Basic properties} \label{sec:basic}
Let $N\in \N$ be fixed. The following properties hold:

\begin{itemize}
\item $Z^N=(Z^N_t)_{t\geq 0}$ is a bounded martingale.\\

\item $Z^N=(Z^N_t)_{t\geq 0}$ is a bounded martingale.\\
Indeed $\E[Z_{t+1}^N|\mathcal{F}_t]=\E[\frac{t+m}{t+m+1}Z^N_t+\frac{1}{N(t+m+1)}\sum_{i=1}^N Y_t(i)|\mathcal{F}_t]=\frac{t+m}{t+m+1}Z^N_t+\frac{1}{t+m+1}\frac{1}{N}\sum_{i=1}^N [\alpha Z^N_t+(1-\alpha)Z_t(i)]=\frac{t+m}{t+m+1}Z^N_t+\frac{1}{t+m+1}Z^N_t=Z^N_t$\\
As a consequence,  $Z^N$ converges a.s. and in $L^p$ to a random variable which we denote by
$Z^N_\infty$.

\item By symmetry reasons, for every fixed $t$, the random variables $Z_t^N(i)$ for $i=1,\ldots,N$ are exchangeable,
thus their mean value equals $\frac{a}{m}$.\\ Indeed, since $Z^N$ is a martingale, for every $t$ we have
  $\frac{a}{m}=\E[Z_t^N]=\E[\frac{1}{N} \sum_{i=1}^N Z_t^N(i)]=\frac{1}{N} \sum_{i=1}^N\E[Z_t^N(i)]=\E[Z_t^N(1)]$.\\
Note that we have also $\E[Y_t(i)]=\E[\alpha Z_t^N+(1-\alpha)Z_t^N(i)]=\frac{a}{m}$ for all
$t,i\in \N$

\item $Z^N(i)=(Z_t^N(i))_{t\geq 0}$ is a martingale if and only if  $\alpha=0$.\\
Indeed $\E[Z_{t+1}^N(i)|\mathcal{F}_t]=\frac{t+m}{t+m+1} Z_t^N(i)+\frac{1}{t+m+1}[\alpha Z_t^N+(1-\alpha)Z_t^N(i)]=
\frac{t+m+1-\alpha}{t+m+1}Z_t^N(i)+\frac{\alpha}{t+m+1}Z_t^N=Z_t^N(i) $ for every $t$ if and only if $\alpha=0$ or $Z_t^N(i)=Z_t^N\ \ \mathbb{P}$-a.s. for every $t$;\ but the last assertion is false, since, for example, when $t=1$ the random variables $Y_0(k)$ for $k=1,\ldots, N$ are i.i.d Bernoulli of parameter $a/m$ and so $\mathbb{P}(Z_1^N(i)=Z_1^N)=\mathbb{P}\big(Y_0(i)=\frac{1}{N}\sum_{k=1}^N Y_0(k)\big)<1$.

\end{itemize}

In next section we study the behavior of $Z^N_t(i)$ for $i=1,\ldots,N$ when $t\rightarrow +\infty$. Since $N$ is fixed, we omit the superscript $N$ and write $Z_t$ and $Z_t(i)$  for $Z^N_t$ and $Z^N_t(i)$ respectively.

\section{Synchronization} \label{sec:sync}

In this section we show that, as soon as $\a>0$, the urns {\em synchronize almost surely}, i.e. for each $i=1,\ldots,N$,
\[
\lim_{t \ra +\infty} Z_t(i) = \lim_{t \ra +\infty} Z_t
\]
almost surely. We proceeds in two steps. First we show $L^2$-synchronization, i.e.
\[
\lim_{t \ra +\infty} \E\left[\left(Z_t(i) - Z_t \right)^2\right] = 0.
\]
Then, using  bounds on the $L^2$-rate of convergence, we derive the almost sure synchronization.

\subsection{$L^2$-synchronization}

In the following statement, for two positive sequences $a_t$, $b_t$ we write $a_t \sim b_t$ if
\[
0 < \liminf_{t \ra +\infty} \frac{a_t}{b_t} \leq \limsup_{t \ra +\infty} \frac{a_t}{b_t} < +\infty.
\]
\bt{th-L^2-sync}
The following asymptotic estimates hold:
\[
\E\left[\left(Z_t(i) - Z_t \right)^2\right] \sim \left\{
\begin{array}{ll}
t^{-2\a} & \mbox{for } 0<\a<\frac12 \\
t^{-1} \log t & \mbox{for } \a = \frac12 \\
t^{-1} & \mbox{for } \frac12 < \a \leq 1.
\end{array}
\right.
\]
\et
For the proof of Theorem 3.1, a technical Lemma is needed.
\bl{lemma:zerovar}
The exists a constant $C$ with $0<C<a/m^2$ such that, for every $t \geq 0$,
\[
\Var\left(Z_t  \right) \leq \frac{C}{N}.
\]
In particular we have $$ \sup_t E(Z_t^2) = E(Z_{\infty}^2) < a/m$$
\el
\begin{proof}[Proof of Lemma~\ref{lemma:zerovar}]
Note that
\[
\mathrm{Var}[Z_{t+1}]=\Var[\E(Z_{t+1}|\mathcal{F}_t)]+\E[\Var(Z_{t+1}|\mathcal{F}_t)] =  \Var[Z_t]+\E[\Var(Z_{t+1}|\mathcal{F}_t)].
\]
Moreover
\begin{eqnarray*}
&&\E[\Var(Z_{t+1}|\mathcal{F}_t)]=
\E\Big[\Var\Big(\frac{t+m}{t+m+1}Z_t+\frac{1}{(t+m+1)N}\sum_{i=1}^N Y_t(i)\Big|\mathcal{F}_t \Big)\Big]\\
&&=\E\Big[\Var\Big(\frac{1}{(t+m+1)N}\sum_{i=1}^N Y_t(i)\Big|\mathcal{F}_t \Big)\Big]\\
&&=\E\Big[\frac{1}{(t+m+1)^2 N^2}\sum_{i=1}^N\Var(Y_t(i)|\mathcal{F}_t)\Big]\\
&&=\E\Big[\frac{1}{(t+m+1)^2 N^2}\sum_{i=1}^N \big(\alpha Z_t+(1-\alpha)Z_t(i)\big) -\big(\alpha Z_t+(1-\alpha)Z_t(i)\big)^2 \Big]\\
&&=\frac{1}{(t+m+1)^2 N}\E\Big[Z_t-\alpha^2 Z_t^2-2(\alpha-\alpha^2)Z_t^2-\frac{(1-\alpha)^2}{N}\sum_{i=1}^N Z_t(i)^2\Big]\\
\end{eqnarray*}
Thus, the following recursive equation holds
\be{eq:riccor}
\begin{split}
\mathrm{Var}[Z_{t+1}]=& \Var[Z_t]+\frac{1}{(t+m+1)^2N}\frac{a}{m}-\frac{\alpha (2-\alpha )}{(t+m+1)^2 N}\Big(\Var[Z_t]+\frac{a^2}{m^2}\Big)\\
& -\frac{(1-\alpha)^2}{(t+m+1)^2 N^2}\sum_{i=1}^N \E[Z_t(i)^2]
\end{split}
\ee
In particular,
\[
\mathrm{Var}[Z_{t+1}]  \leq  \Var[Z_t]+\frac{1}{N(t+m+1)^2}\frac{a}{m},
\]
from which we get
$$
\Var[Z_t] \leq \frac{a}{m}\frac{1}{N} \sum_{k=m+1}^{+\infty} \frac{1}{k^2} < \frac{a}{m^2N}.
$$
Moreover
\[
\E(Z_t^2) = \frac{a^2}{m^2} + \Var[Z_t] <  \frac{a^2 +a}{m^2} \leq \frac{a}{m},
\]
where this last inequality is easily shown to be equivalent to $a \leq m-1$, which is clearly true.
Finally, since $Z_t^2$ is a bounded submartingale,
\[
 \sup_t \E(Z_t^2) = \lim_t \E(Z_t^2) = \E(Z_{\infty}^2).
\]
\epr
\begin{proof}[Proof of Theorem \ref{th-L^2-sync}]
Setting
\[
x_t := \E\left[\left(Z_t(i) - Z_t \right)^2\right] = \Var(Z_t(i) - Z_t),
\]
we first obtain a recursive equation satisfied by $x_t$.
\begin{multline*}
\begin{split}
x_{t+1}&=\E\big[\Var\big(Z_{t+1}(i)-Z_{t+1}|\mathcal{F}_t\big)\big]+\Var\big(\E\big[Z_{t+1}(i)-Z_{t+1}|\mathcal{F}_t\big] \big)\\
&=\E\left[\Var\left(\frac{1}{t+m+1}\Big\{(t+m)\big(Z_t(i)-Z_t\big)+Y_{t}(i)-\frac{1}{N}\sum_{i=1}^N Y_t(i)\Big\}\Big|\mathcal{F}_t\right)\right]\\
& \ \ \ +\Var\left( \frac{t+m}{t+m+1}Z_t(i)+\frac{1}{t+m+1}[\alpha Z_t+(1-\alpha)Z_t(i)]-Z_t\right)\\
&= \frac{1}{(t+m+1)^2}\E\Big[\Var\Big(Y_t(i)-\frac{1}{N}\sum_{j=1}^N Y_t(j)|\mathcal{F}_t\Big)\Big]\\
& \ \ \ + \Var\Big(\frac{t+m+1-\alpha}{t+m+1}\big[Z_t(i)-Z_t\big]\Big)\\
&=\frac{1}{(t+m+1)^2}\Big[\Big(1-\frac{1}{N}\Big)^2+\frac{N-1}{N^2}\Big]\E\big[\Var\big(Y_t(i)|\mathcal{F}_t \big)\big]\\
& \ \ \ +\Big(\frac{t+m+1-\alpha}{t+m+1}\Big)^2\Var\big(Z_t(i)-Z_t(i)\big)\\
&=\frac{1}{(t+m+1)^2}\frac{N-1}{N}\Big(\frac{a}{m}-\E\big[\big(\alpha Z_t+(1-\alpha)Z_t(i)\big)^2\big]\Big)\\
& \ \ \ +\Big(1-\frac{\alpha}{t+m+1}\Big)^2x_t\\
&=\frac{1}{(t+m+1)^2}\frac{N-1}{N}\Big(\frac{a}{m}-\E\big[\big\{Z_t-(1-\alpha)\big(Z_t-Z_t(i)\big)\big\}^2\big]\Big)\\
& \ \ \ +\Big(1-\frac{\alpha}{t+m+1}\Big)^2x_t\\
&=\frac{1}{(t+m+1)^2}\frac{N-1}{N}\Big(\frac{a}{m}-(1-\alpha)^2 x_t-\E\big[Z_t^2\big]\Big)\\
& \ \ \ +x_t-\frac{2\alpha}{(t+m+1)}x_t+\frac{\alpha^2}{(t+m+1)^2}x_t
\end{split}
\end{multline*}
where, in the last equality, we have used the fact that $\E\big[Z_t\big(Z_t-Z_t(i)\big)\big]=0$ since, by symmetry, $$\E\big[Z_t(i)Z_t\big]=\frac{1}{N}\sum_{j=1}^N\E\big[Z_t(j)Z_t\big]=
\E\big[\big(\frac{1}{N}\sum_{j=1}^N Z_t(j)\big)Z_t\big]=\E[Z_t^2].$$
We have therefore obtained the recursive equation
\be{varianza}
\begin{split}
x_{t+1}= & \ x_t-\frac{2\alpha}{t+m+1}x_t+
\frac{\alpha^2-\frac{N-1}{N}(1-\alpha)^2}{(t+m+1)^2}x_t
\\ & +\frac{1}{(t+m+1)^2}\frac{N-1}{N}\Big(\frac{a}{m}-\E[Z_t^2]\Big)\\
\end{split}
\ee
Now set $A := 2\a$, $B := \alpha^2-\frac{N-1}{N}(1-\alpha)^2$,
\[
f(t) := 1- \frac{A}{t+m+1} + \frac{B}{(t+m+1)^2}, \ \ \ g(t):= \frac{\frac{N-1}{N}\Big(\frac{a}{m}-\E[Z_t^2]\Big)}{(t+m+1)^2},
\]
so that
\be{x-rec}
x_{t+1} = f(t) x_t + g(t).
\ee
It is easily seen that $-1 \leq B \leq 1$ and, since the previous lemma~\ref{lemma:zerovar} states $\frac{a}{m}-\E[Z_{\infty}^2]>0$, we have
\begin{equation}\label{GT}
0< \frac{\frac{N-1}{N}\Big(\frac{a}{m}-\E[Z_{\infty}^2]\Big)}{(t+m+1)^2} \leq g(t) \leq \frac{\frac{a}{m}}{(t+m+1)^2}. \
\end{equation}

We also remark that $0 < f(t) < 1$ for every $t \geq 0$. To see this, observe first that
\[
f(0) = 1-\frac{2\a}{m+1} + \frac{\a^2 - \frac{N-1}{N}(1-\a)^2}{(m+1)^2}
\]
as a function of $\a$, has a strictly negative derivative; being positive for $\a=1$, it must be positive for all $\a \in [0,1]$. Moreover
\begin{multline*}
\begin{split}
\frac{d}{dt} f(t)&= \frac{2}{(t+m+1)^3} \left[ \a(t+m+1) - \a^2 + \frac{N-1}{N} (1-\a)^2 \right]\\& > \frac{2}{(t+m+1)^3} (\a - \a^2) \geq 0,
\end{split}
\end{multline*}
so $f(t)$ is increasing in $t$. Since $\lim_{t \ra +\infty} f(t) = 1$, we conclude $0 < f(t) < 1$. Now set
\[
\xi_t := \frac{x_t}{\prod_{k=0}^{t-1} f(k)}.
\]
By \eqref{x-rec}, we obtain
\[
\xi_{t+1} = \xi_t + F(t),
\]
where
\[
F(t) := \frac{g(t)}{\prod_{k=0}^{t} f(k)}.
\]
So, observing that $\xi_0 = x_0 = 0$, we get
\[
\xi_t = \sum_{i=0}^{t-1} F(i)
\]
or, equivalently,
\be{x_t}
x_t = \left[\prod_{k=0}^{t-1} f(k)\right]  \sum_{i=0}^{t-1} F(i).
\ee
Now we can study the asymptotic behavior of each term. Note first that, as $t \ra +\infty$,
\be{x_t-1}
\begin{split}
 \prod_{k=0}^{t-1} f(k) & = \exp\left[ \sum_{k=0}^{t-1} \log \left(1-\frac{A}{k+m+1} + \frac{B}{(k+m+1)^2} \right) \right]\\
& = \exp\left[ -2\a  \sum_{k=0}^{t-1} \frac{1}{k+m+1} + O(1) \right]\\
& = \exp[-2\a \log (t+m) + O(1)] \sim t^{-2\a}.
 \end{split}
 \ee
Since, by (\ref{GT}), we have $g(t) \sim t^{-2}$,
\be{x_t-2}
F(t) \sim t^{2\a - 2} \ \Rightarrow \ \xi_t  = \sum_{i=0}^{t-1} F(i) \sim \left\{ \begin{array}{ll}  1 & \mbox{for } 0 \leq \a < \frac{1}{2} \\ \log t & \mbox{for } \a = \frac12 \\ t^{2\a -1} &  \mbox{for } \frac12 < \a \leq 1. \end{array} \right.
\ee
By \eqref{x_t}, \eqref{x_t-1} and \eqref{x_t-2} the conclusion follows.

\end{proof}

\subsection{Almost sure synchronization}
For $\a = 0$, the sequences $\{Z_t(i)\}_{t\geq 0}$ converge, as $t \ra +\infty$, to independent and Beta distributed limits $Z_{\infty}(i)$. For $\a>0$ we have a quite different behavior.
\bt{th-a.s.-sync}
For each $i=1,2,\ldots,N$,
\[
\lim_{t \ra +\infty} Z_t(i) = \lim_{t \ra +\infty} Z_t
\]
almost surely.
\et

\bpr
By Theorem \ref{th-L^2-sync}, $Z_t(i)$ converges to $Z_{\infty}$ in $L^2$, as $t\mapsto \infty$. It is therefore enough to show that the almost sure limit $\lim_{t \ra +\infty} Z_t(i)$ exists. We observe that
$\{Z_t(i)\}_t$ is a quasi-martingale, i.e.
\[
 \sum_{t=0}^{+\infty} \E \left[\ \left| \E(Z_{t+1}(i) | \mathcal F_t) - Z_t(i)\right|\ \right] < + \infty .
 \]
 Indeed, using the fact that
 \[
\E(Z_{t+1}(i) | \mathcal F_t) = \left(1-  \frac{\a}{t+m+1} \right) Z_t(i) + \frac{\a}{t+m+1} Z_t,
\]
we obtain
\[
 \begin{split}
 \sum_t \E \left[\ \left| \E(Z_{t+1}(i) | \mathcal F_t) - Z_t(i)\right|\  \right]  & =  \sum_t \frac{\a}{t+m+1}\E  \left[ \ \left| Z_t(i) - Z_t\right| \right] \\ & \leq \sum_t \frac{\a}{t+m+1} \left(\E\left[\left(Z_t(i) - Z_t \right)^2\right]\ \right)^{1/2} < +\infty,
 \end{split}
 \]
 where the last inequality comes from Theorem \ref{th-L^2-sync}. The conclusion now follows from the fact that a bounded quasi-martingale has an almost sure limit (see e.g. \cite{Metivier} pag. 46)

\epr

\section{A central limit theorem} \label{sec:CLT}

In this section we study the limiting distribution of $Z_t$ in the limit as the number $N$ of urns goes to infinity. To emphasize the dependence on the number of urns, we come back to notation $Z_t^{N}$ and set
\[
W_t^{N} := \sqrt{N} \left( Z^{N}_t - \frac{a}{m}\right).
\]
Moreover, let $x_t^{\infty}$ be the solution of the recursion
\begin{eqnarray} \label{eq:x-infty}
x^\infty_{t+1}=x^\infty_t-\frac{2\alpha}{t+m+1}x^\infty_t + \frac{(2\alpha-1 )}{(t+m+1)^2} x^\infty_t %
+\frac{1}{(t+m+1)^2} \Big(\frac{a}{m}-\frac{a^2}{m^2}\Big),
\end{eqnarray}
with $x^{\infty}_0 = 0$. Note that \eqref{eq:x-infty}  is the limit as $N \rightarrow +\infty$ of \eqref{varianza}, which implies
\[
x^\infty_t = \lim_{N \rightarrow +\infty} \E\left[ (Z_t^N(i) - Z_t^N)^2 \right].
\]
\bt{th:clt}
The stochastic process $\{W_t^{N} : \ t \geq 0\}$ converges weakly, as \mbox{$N \ra +\infty$},
to the Gau\ss-Markov process solution of the recursion
\begin{equation} \label{eq:CLT}
\left\lbrace \begin{array}{l}
 W_{t+1} = W_{t} +  \sigma_t B_{t+1} \\
W_0= \delta_0
\end{array} \right.
\end{equation}
where
\be{sigmat}
\sigma_t^2=\frac{1}{(t+m+1)^2} \left(\left[\frac{a}{m}-\frac{a^2}{m^2}\right] - (1-\alpha)^2 x^{\infty}_t \right) \geq 0,
\ee
 and $\{B_{t}: \ t \geq 1\}$ is a sequence of i.i.d. $\mathcal N(0,1)$.

\et
We begin with two technical lemmas.
\bl{lemma:lln}
The following almost sure limits hold, for every $t \geq 1$:
\be{clt1}
\lim_{N \ra +\infty} Z_t^{N} = \frac{a}{m},
\ee
\be{clt2}
\lim_{N \ra +\infty} \frac{1}{N} \sum_{i=1}^N \left[Z_t^{N}(i) \right]^2 = x_t^{\infty} + \frac{a^2}{m^2}.
\ee

\el

\bpr
We begin by proving \eqref{clt1}, by induction on $t$. Since $Z_0^{N} \equiv \frac{a}{m}$, there is nothing to prove for $t=0$. For simplicity, in next formulas we go back to the notations $Z_t$, $Z_t(i)$, omitting the dependence on $N$. Recall that
\be{clt3}
Z_{t+1} = Z_t + \frac{1}{t+m+1} \frac{1}{N} \sum_{i=1}^N \left[ Y_t(i) - \a Z_t - (1-\a) Z_t(i) \right].
\ee
Let
\[
V_i :=  Y_t(i) - \a Z_t - (1-\a) Z_t(i) .
\]
By assumption the $V_i$'s are independent for the conditional probability $\Prob\left(\, \cdot \, | \mathcal{F}_t \right)$ and $\E\left(V_i |\mathcal{F}_t\right) =0$. Since the $V_i$'s are bounded, a strong law of large numbers holds: letting
\[
F:= \left\{ \lim_{N \ra +\infty} \frac{1}{N} \sum_{i=1}^N V_i = 0 \right\},
\]
we have $\Prob\left(F | \mathcal{F}_t \right) = 1$ a.s., which implies $\Prob(F) = 1$. Thus, by \eqref{clt3},
\[
\lim_{N \ra +\infty} Z_{t+1} = \lim_{N \ra +\infty} Z_t = \frac{a}{m}
\]
almost surely, where we have used the inductive assumption.

We now prove \eqref{clt2}, again by induction on $t$. Using
\[
Z_{t+1}(i) = \frac{t+m}{t+m+1} Z_t(i) + \frac{1}{t+m+1} Y_t(i),
\]
we get
\[
\begin{split}
Z^2_{t+1}(i)  = & \left( \frac{t+m}{t+m+1}\right)^2 Z^2_t(i) + \frac{1}{(t+m+1)^2} Y_t(i) + \frac{2(t+m)}{(t+m+1)^2}Z_t(i) Y_t(i) \\
= & \left( \frac{t+m}{t+m+1}\right)^2 Z^2_t(i) + \frac{1}{(t+m+1)^2} V_i +
\frac{\alpha Z_t + (1-\alpha)Z_t(i)}{(t+m+1)^2} \\
 & + \frac{2(t+m)}{(t+m+1)^2} Z_t(i) V_i + \frac{2(t+m)}{(t+m+1)^2}Z_t(i)\big[\alpha Z_t+(1-\alpha)Z_t(i)\big],
\end{split}
\]
which gives
\be{clt4}
\begin{split}
&\frac{1}{N} \sum_{i=1}^N Z^2_{t+1}(i)  = \frac{(t+m)(t+m+2-2\alpha)}{(t+m+1)^2} \frac{1}{N} \sum_{i=1}^N Z^2_{t}(i) + \frac{1}{(t+m+1)^2} Z_t \\
& + \frac{2 \a (t+m)}{(t+m+1)^2}  Z^2_t +\frac{1}{(t+m+1)^2}\frac{1}{N} \sum_{i=1}^N  V_i + \frac{2(t+m)}{(t+m+1)^2}\frac{1}{N} \sum_{i=1}^N Z_t(i) V_i .
\end{split}
\ee
The fact that $\frac{1}{N} \sum_{i=1}^N Z^2_{t+1}(i)$ converges almost surely, as $N$ tends to infinity, to a constant follows from \eqref{clt4}, by:
\begin{itemize}
\item
the inductive assumption;
\item
the fact, already proved, that {$\lim_{N\to \infty} Z_t = \frac{a}{m}$} a.s.; 
\item
the fact that
\[
\frac{1}{N} \sum_{i=1}^N  V_i  \ \ra \ 0 \ \mbox{a.s.}
\]
as shown in the proof of \eqref{clt1};
\item
the fact that
\[
\frac{1}{N} \sum_{i=1}^N Z_t(i) V_i \ \ra \ 0 \ \mbox{a.s.},
\]
which is proved in the same way.
\end{itemize}
To identify this limit, note that
\[
\frac{1}{N}\sum_{i=1}^N\E[Z_t(i)^2]=\Var[Z_t-Z_t(i)]+\E[Z_t^2]=x_t+\Var[Z_t]+\frac{a^2}{m^2},
\]
By Lemma~\ref{lemma:zerovar}, $\Var[Z_t]$ tends to zero as $N \ra +\infty$. Finally, since $x_t$ satisfies \eqref{varianza}, it is easily proved by induction on $t$ that $x_t$ has a limit $x_t^{\infty}$ as $N \ra +\infty$, and this limit satisfies \eqref{eq:x-infty}.

\epr

\begin{remark} \label{rem:sigmat}
From Lemma \ref{lemma:lln} it follows that $\sigma_t^2 \geq 0$, where $\sigma_t$ has been defined in \eqref{sigmat}. Indeed
\[
\frac{a}{m} = \lim_{N \ra +\infty} \frac{1}{N} \sum_{i=1}^N Z_t^{N}(i) \geq \lim_{N \ra +\infty} \frac{1}{N} \sum_{i=1}^N \left[Z_t^{N}(i) \right]^2 = x_t^{\infty} + \frac{a^2}{m^2},
\]
thus
\[
x_t^{\infty} \leq \frac{a}{m} - \frac{a^2}{m^2},
\]
which implies $\sigma_t^2 \geq 0$.

\end{remark}

\begin{lemma} \label{lemma:chf}
$$ \E (e^{\textrm{i} u W_{t+1}^{N}} | \mathcal F_t) = e^{\textrm{i} u W_{t}^{N}} %
e^{-\frac{u^2}{2}\frac{1}{(t+m+1)^2}\left[ (\frac{a}{m}-\frac{a^2}{m^2})-(1-\alpha)^2x_t^\infty + \varepsilon_N \right]}$$
where $\varepsilon_N$ is a sequence of $\mathcal{F}_t$-measurable, uniformly bounded random variables, such that $\varepsilon_N \rightarrow 0$ a.s.
\end{lemma}
\begin{proof} Using the recursive equation
\[
W_{t+1}^{N} = \frac{t+m}{t+m+1}W_t^{N} + \frac{1}{t+m+1} \left(\frac{1}{\sqrt{N}} \sum_{i=1}^N Y_t(i)\ -\frac{a}{m}\sqrt{N} \right),
\]
and using the notation $s=m+t$ we have
$$ \E (\textrm{e}^{\textrm{i} u W_{t+1}^{N}} | \mathcal F_t) = \textrm{e}^{i u \frac{s}{s+1}W_{t}^{N}} %
\textrm{e}^{-\textrm{i}  \frac{u}{s+1} \frac{a\sqrt{N}}{m}} \  \E  \left(\prod_{j=1}^N  \textrm{e}^{\textrm{i} \frac{u}{\sqrt{N} (s+1)} Y_t(j)} \Big| \mathcal{F}_t \right) .$$
Compute, using the characteristic function of a Bernoulli distribution,
{\footnotesize{
\begin{multline*}
\begin{split}
  &\E  \left(\prod_{j=1}^N  \textrm{e}^{\textrm{i} \frac{u}{\sqrt{N} (s+1)} Y_t(j)} | \mathcal{F}_t\right) =  %
 \prod_{j=1}^N \big[  1+ (\textrm{e}^{\textrm{i} \frac{u}{\sqrt{N} (s+1)}} -1) (\alpha Z_t + (1-\alpha ) Z_t(j)) \big] \\
 =& \exp \left( \sum_{j=1}^N \log \left[ 1 + \left( \frac{\textrm{i} u}{\sqrt{N} (s+1)} %
- \frac{u^2}{2 N (s+1)^2} + o(\frac{1}{N}) \right) \ (\alpha Z_t + (1-\alpha) Z_t(j) )\right]   \right) \\
 =& \exp \left(  \frac{\textrm{i} u \sqrt{N}}{s+1} Z_t -  \frac{u^2}{2  (s+1)^2} Z_t + \frac{u^2}{2(s+1)^2} \frac{1}{N} \sum_{j=1}^N %
\Big( \alpha Z_t + (1-\alpha) Z_t(j) \Big)^2 + o(1) \right),
\end{split}
\end{multline*}
}}
where, in this last expression, $o(1)$ denotes a uniformly bounded sequence of random variables which goes to zero a.s. as $N \rightarrow +\infty$ that, in what follows, may change from line to line.\\
Recalling  $Z_t= \frac{W_t^{(N)}}{\sqrt{N}} + \frac{a}{m}$, the first term in the argument of the exponential above equals $iu\frac{1}{s+1}W_t+\textrm{i}\frac{u}{s+1}\frac{a\sqrt{N}}{m}$. Thus
{\footnotesize{
\begin{multline*}
\E (\textrm{e}^{\textrm{i} u W_{t+1}^{N}} | \mathcal F_t)  \\ = \textrm{e}^{\textrm{i} u W_{t}^{N}} \exp \left(
-\frac{u^2}{2(s+1)^2} \Big[\ Z_t-\alpha^2 Z_t^2-2\alpha(1-\alpha)Z_t^2-(1-\alpha)^2\frac{1}{N}\sum_{j=1}^N Z^2_t(j) \ \Big] + o(1) \right).
\end{multline*}
}}
which, using  Lemma \ref{lemma:lln}, yields:
$$ \E (\textrm{e}^{\textrm{i} u W_{t+1}^{N}} | \mathcal F_t) = \textrm{e}^{\textrm{i} u W_{t}^{N}} \exp \left(-\frac{u^2}{2}
\frac{1}{(s+1)^2} \Big[\frac{a}{m}-\frac{a^2}{m^2}-(1-\alpha)^2x_t^\infty \Big] +o(1) \right) .$$

\end{proof}

\begin{proof}[Proof of Theorem \ref{th:clt}]
Since $$\E\Big[\exp \Big( \textrm{i} \sum_{r=1}^{t+1}  u_r W^{N}_r \Big)\Big] = \E\Big[\exp \Big( \textrm{i} \sum_{r=1}^t  u_r W^{N}_r \Big) \ \E \Big(\exp \big( \textrm{i}\ u_{t+1}  W^{N}_{t+1} \big) \Big| \mathcal F_t \Big)\Big]$$
by Lemma \ref{lemma:chf} we easily see, by induction on $t$, that the limits
\[
\varphi_t(u_1,\cdots,u_t) := \lim_{N \to \infty} \E\Big[\exp \Big( \textrm{i} \sum_{r=1}^t u_r W^{N}_r \Big)\Big)
\]
exist, and satisfy the recursion
\begin{multline*}
 \varphi_{t+1} (u_1,\ldots,u_{t-1},u_t,u_{t+1})  \\ = \varphi_{t} (u_1,\ldots,u_{t-1},u_t+u_{t+1}) \ \exp \left(-\frac{u^2}{2{(t+m+1)^2}} \Big[\frac{a}{m}-\frac{a^2}{m^2}-(1-\alpha)^2x_t^\infty \Big]\right),
 \end{multline*}
with initial condition $\varphi_0 =1$. It is now straightforward to check that the Gaussian process in \eqref{eq:CLT} gives rise to the same recursion for the characteristic functions.

\end{proof}

\br{Markov}
Note that, for $N$ finite, $(W_t^{N})$ is not a Markov process: the whole vector $(Z_t(i))_{i=1}^N$ is needed to have a Markov process. The Markov property is recovered in the limit as $N \rightarrow +\infty$.
\er


\end{document}